\documentclass[english]{article}
\usepackage[T1]{fontenc}
\usepackage[latin9]{inputenc}
\usepackage{textcomp}
\usepackage{relsize}
\usepackage{amsmath}
\usepackage{amssymb}

\makeatletter

\DeclareRobustCommand{\lyxmathsym}[1]{\ifmmode\begingroup\def\b@ld{bold}
  \def\rmorbf##1{\ifx\math@version\b@ld\textbf{##1}\else\textrm{##1}\fi}
  \mathchoice{\hbox{\rmorbf{#1}}}{\hbox{\rmorbf{#1}}}
  {\hbox{\smaller[2]\rmorbf{#1}}}{\hbox{\smaller[3]\rmorbf{#1}}}
  \endgroup\else#1\fi}

\providecommand{\tabularnewline}{\\}

\makeatother

\makeatother

\usepackage{babel}

\makeatother

\usepackage{babel}

\makeatother

\usepackage{babel}

\makeatother

\usepackage{babel}

\makeatother

\usepackage{babel}

\makeatother

\usepackage{babel}

\makeatother

\usepackage{babel}

\makeatother

\usepackage{babel}

\makeatother

\usepackage{babel}

\makeatother

\usepackage{babel}

\makeatother

\usepackage{babel}

\makeatother

\usepackage{babel}

\begin{document}
\numberwithin{equation}{section}

\title{On the Distribution of the non-trivial Roots of Riemann's Zeta Function.
An alternative Approach.}

\author{Renaat Van Malderen}
\maketitle
\begin{abstract}
\noindent Expressing Weierstrass type infinite products in terms of
Stieltjes integrals is discussed. The asymptotic behavior of particular
types of infinite products is compared against the asymptotic behavior
of the entire function $\xi$(s), well-known in Riemann zeta function
theory. An approximate formula for the distribution of the non-trivial
roots of Riemann's zeta function is obtained.\\

\end{abstract}

\section{\textit{Introduction.}}

\noindent In the theory of Riemann's zeta function the entire function

\noindent \begin{equation}
\xi\left({\rm s}\right){\rm =}\frac{1}{2}\Gamma\left(\frac{{\rm s}}{{\rm 2}}\right){\pi}^{{\rm -}\frac{{\rm s}}{{\rm 2}}}{\rm s}\left({\rm s-1}\right)\zeta\left({\rm s}\right)\label{eq:section_1_1}\end{equation}
 plays a fundamental role. As demonstrated in this theory, the zeros
of $\xi$(s) occur only in the critical strip $0<\Re(s)<1$, and given
the so called Riemann hypothesis is true, these zeros will all be
located on the vertical axis $\Re(s)=\frac{1}{2}$.

\noindent Theorems developed by Hadamard in 1893 allow to demonstrate
that $\xi(s)$ may be represented by an infinite Weierstrass product
(see e.g. Lang, p. 386):

\noindent \begin{equation}
\xi{\rm (s)=}\xi{\rm (0)}\prod_{\rho}{\left({\rm 1-}\frac{{\rm s}}{\rho}\right)}\label{eq:section_1_2}\end{equation}

\noindent in which the set $\rho$ represents the zeros of $\xi(s)$,
which are all located within the critical strip (see e.g. Edwards,
p.18). An important point to note in the product \eqref{eq:section_1_2}
is that no exponential convergence factor is required per individual
term. $\xi(s)$ is analytic in the entire $s$ plane and is an even
function of $s-\frac{1}{2}$, i.e. $\xi(s-\frac{1}{2})=\xi(\frac{1}{2}-s)$,
and its roots occur in pairs, i.e. if the complex number $\rho$ is
a root, so is $1-\rho$.

\noindent In \eqref{eq:section_1_2} the terms in the product are
to be taken in pairs, i.e. $\rho$ together with $1-\rho$.

\section{\textit{Types of infinite products.}}

\noindent In what follows we will consider two types of infinite products
which are a generalisation of \eqref{eq:section_1_2}: 
\begin{enumerate}
\item Pure axial products in which all zeros lie on the axis $\Re(s)=\frac{1}{2}$. 
\item Strip type products in which some zeros also lie outside the same
axis as in 1 but still within the critical strip. 
\end{enumerate}
\noindent The axial products are a subset of the strip product set.
$\xi(s)$ belongs at least to the strip product set and if the Riemann
hypothesis were true, $\xi(s)$ would be an axial product. For both
1 and 2 we exclude roots at $s=\frac{1}{2}$.

\noindent Below, instead of the complex variable $s=\sigma+it$ usual
in analytic number theory, we will use:

\[
{\rm z=s-}\frac{{\rm 1}}{{\rm 2}}=x+it\]
 This brings more symmetry into the considerations. For both axial
and strip types we will consider products F(z) analytic in the entire
z plane (entire functions) which are in addition even, i.e.

$F(z)=F(-z)$

\noindent which implies $F(x)=F(-x)$ and $F(it)=F(-it)$.

\noindent We also assume $F(x)$ as well as $F(it)$ to be real and
the product expansion for $F(z)$ requiring no exponential convergence
factors. We will therefore consider: 
\begin{enumerate}
\item Axial products: their zeros are purely imaginary and occurring in
pairs $\pm z_{l}$with: \begin{equation}
z_{l}=ik_{l}\quad k_{l}>0\label{eq:section_2_1}\end{equation}
 The product\begin{equation}
{\rm F}\left({\rm z}\right){\rm =F}\left(0\right)\prod_{{\rm l=1}}^{\infty}{{\rm (1+\ }\frac{{\rm z}^{{\rm 2}}}{{\rm k}_{{\rm l}}^{{\rm 2}}})}\label{eq:section_2_2}\end{equation}
 with F(0) real, represents an even entire function in the z plane.
As well known \eqref{eq:section_2_2} converges if\begin{equation}
\sum_{{\rm l=1}}^{\infty}{\frac{{\rm 1}}{{\rm k}_{{\rm l}}^{{\rm 2}}}}\label{eq:section_2_3}\end{equation}
 converges 
\item Strip products: their off-axis zeros need to occur in quads, i.e.
$z_{m}$, $-z_{m}$, $\overline{z_{m}}$,-$\overline{z_{m}}$, all
four values lying within the critical strip. We will use the representations: 
\end{enumerate}
\noindent \[
z_{m}=iq_{m}e^{-i\beta m}\]

\begin{equation}
-z_{m}=-iq_{m}e^{-i\beta m}\label{eq:section_2_4}\end{equation}

\[
\overline{z_{m}}=-iq_{m}e^{i\beta m}\]

\[
\overline{{-z}_{m}}=iq_{m}e^{i\beta m}\]
 in which we assume $q_{m}$ to be real and $>\frac{1}{2}$

The associated product then equals:

\noindent \[
{\rm G}\left({\rm z}\right){\rm =G}\left(0\right)\prod_{{\rm m=1}}^{\infty}{\left({\rm 1-}\frac{{\rm z}}{{\rm z}_{{\rm m}}}\right)}\left({\rm 1+}\frac{{\rm z}}{{\rm z}_{{\rm m}}}\right)\left({\rm 1-}\frac{{\rm z}}{{\overline{{\rm z}_{{\rm m}}}}}\right)\left({\rm 1+}\frac{{\rm z}}{{\overline{{\rm z}_{{\rm m}}}}}\right)\]

\begin{equation}
{\rm G}\left({\rm z}\right){\rm =G}\left(0\right)\prod_{{\rm m=1}}^{\infty}{\left({\rm 1+\ }\frac{{\rm 2}{\rm z}^{{\rm 2}}{\rm cos}{\beta}_{{\rm m}}}{{\rm q}_{{\rm m}}^{{\rm 2}}}{\rm +}\frac{{\rm z}^{{\rm 4}}}{{\rm q}_{{\rm m}}^{{\rm 4}}}\right)}\label{eq:section_2_5}\end{equation}
 with $G(0)$ real.

\noindent As for axial products we assume $\sum_{{\rm m=1}}^{\infty}{\frac{{\rm 1}}{{\rm q}^{{\rm 2m}}}}$
to converge. Strip products then consist of the combined product of
expressions \eqref{eq:section_2_2} and \eqref{eq:section_2_5}. In
what follows we will first consider axial products. Later we will
look at strip products and multiple roots.

\section{\textit{Axial products.}}

\noindent By taking the logarithm of \eqref{eq:section_2_2} we split
the product into a sum:

\noindent \begin{equation}
{\rm f}\left({\rm z}\right){\rm =}ln{\rm F}\left({\rm z}\right){\rm =}ln{\rm F}\left(0\right){\rm +\ }\sum_{{\rm l=1}}^{\infty}{ln{\rm (1+}\frac{{\rm z}^{{\rm 2}}}{{\rm k}_{{\rm l}}^{{\rm 2}}}}{\rm )}\label{eq:section_3_1}\end{equation}
 ln $F(z)$ is here defined by:

\noindent \[
{\rm f}\left({\rm z}\right){\rm =}ln{\rm F}\left(0\right){\rm +\ }\int_{0}^{{\rm z}}{\frac{{\rm F}^{{\rm '}}{\rm (w)}}{{\rm F(w)}}}{\rm dw}\]

\noindent in the open semiplane $x>0$, via a path of integration
not crossing the axis $x=0$. This way $f(z)$ is defined unambiguously
and the singularities occurring at $z=\pm ik_{l}$ are avoided. The
fact that $F(0)$ is taken at the point $z=0$, i.e. on the axis $x=0$
is no problem for the definition of $f(z)$. In conjunction with the
distribution of the zeros $z=\pm ik_{l}$ we now introduce a non-decreasing
step function for $k\ge0$:

\noindent \begin{equation}
\varphi\left({\rm k}\right){\rm =\ }\sum_{{\rm k}_{{\rm l}}\le{\rm k}}{\rm 1}\label{eq:section_3_2}\end{equation}

\noindent In words: for $k_{l}\leq k<k_{l}+1$, $\varphi(k)$ equals
the number of roots less than or equal to $k_{l}$. The way we have
defined our infinite products, $\varphi(k)$ always equals zero for
$0\leq k<k_{1}$ and then jumps to the value 1. In order for \eqref{eq:section_3_1}
to converge, the infinite sum $\sum_{l=1}^{\infty}{\frac{1}{k_{l}^{2}}}$
should converge (see \eqref{eq:section_2_3}).

\noindent This implies also that there is no accumulation point for
the $k_{l}$ for any finite value of k. We also ignore multiple roots
for the moment. \eqref{eq:section_3_1} may now be expressed as a
Stieltjes integral:

\noindent \begin{equation}
{\rm f}\left({\rm z}\right){\rm =}ln{\rm F}\left(0\right){\rm +\ }\int_{0}^{\infty}{ln{\rm (1+}\frac{{\rm z}^{{\rm 2}}}{{\rm k}^{{\rm 2}}}}{\rm )d}\varphi\left({\rm k}\right)\label{eq:section_3_3}\end{equation}

\noindent Partial integration yields:

\noindent \[
{\rm f}\left({\rm z}\right){\rm =}ln{\rm F}\left(0\right){\rm +}\varphi\left({\rm k}\right){ln{\left.\left({\rm 1+}\frac{{\rm z}^{{\rm 2}}}{{\rm k}^{{\rm 2}}}\right)\right|}_{0}^{\infty}}{\rm -}\int_{0}^{\infty}{\varphi{\rm (k)d}\left[ln{\rm (1+}\frac{{\rm z}^{{\rm 2}}}{{\rm k}^{{\rm 2}}}{\rm )}\right]}\]

\noindent For the second term on the right, in the lower limit $\varphi\left({\rm k}\right)ln{\left.\left({\rm 1+}\frac{{\rm z}^{{\rm 2}}}{{\rm k}^{{\rm 2}}}\right)\right|}_{0}$
equals zero since $\varphi(k)=0$ for $0\leq k<k_{1}$.

\noindent Consider now the upper limit ($k\rightarrow\infty$) of
the same term. For a given value of z and k$\to\infty$:

\noindent \[
{ln\left({\rm 1+}\frac{{\rm z}^{{\rm 2}}}{{\rm k}^{{\rm 2}}}\right)}\sim\frac{{\rm z}^{{\rm 2}}}{{\rm k}^{{\rm 2}}}\]

\noindent As long as for $k\rightarrow\infty$, $\varphi(k)<C\; k^{\alpha}$
with $C$ a positive constant and $\alpha<2$, the term ${\varphi{\rm (k)ln}\left(1+\frac{z^{2}}{k^{2}}\right)}$
goes to zero as well. This is a constraint we put on $\varphi(k)$
and is in fact the equivalent of \eqref{eq:section_2_3} Under these
assumptions:

\noindent \[
{\rm f}\left({\rm z}\right){\rm =}{ln{\rm F}\left(0\right){\rm -}\ }{\rm \ }\int_{0}^{\infty}{\varphi\left({\rm k}\right){\rm d}\left[ln{\rm (1+\ }\frac{{\rm z}^{{\rm 2}}}{{\rm k}^{{\rm 2}}}\right]}\]

\begin{equation}
{\rm f}\left({\rm z}\right){\rm =}{ln{\rm F}\left(0\right){\rm +}\ }{\rm 2}{\rm z}^{{\rm 2}}{\rm \ }\int_{0}^{\infty}{\frac{\varphi{\rm (k)dk}}{{\rm k(}{\rm k}^{{\rm 2}}{\rm +}{\rm z}^{{\rm 2}}{\rm )}}}\label{eq:section_3_4}\end{equation}

\section{\textit{Axial product example: Cosh z}}

\noindent A straightforward example is the function $F(z)=Cosh\: z$,which
for $z=it$ becomes $F(it)=cos\: t$.

\noindent $Cosh\: z$ is an entire function and its zeros occur in
pairs along the imaginary axis at equidistant points:

\noindent \[
{\rm k}_{{\rm l}}{\rm =\pm(2l-1)}\frac{\pi}{{\rm 2}}\qquad(l=1,2,\dots)\]
 and as well known:

\noindent \begin{equation}
{\rm Coshz=}\prod_{{\rm l=1}}^{\infty}{{\rm (1+}\frac{{\rm 4z}^{{\rm 2}}}{{\pi}^{{\rm 2}}{\rm (2l-1)}^{{\rm 2}}})}\label{eq:section_4_1}\end{equation}

\noindent Also

\begin{equation}
{{\rm Cosh}{\rm z=}\frac{{\rm e}^{{\rm z}}{\rm +}{\rm e}^{{\rm -}{\rm z}}}{{\rm 2}}}\label{eq:section_4_2}\end{equation}

\noindent We will now show how to recover \eqref{eq:section_4_2}
from \eqref{eq:section_4_1}.through the use of \eqref{eq:section_3_4}.
In \eqref{eq:section_3_4} since F(0) in \eqref{eq:section_4_1} equals
1, ln F(0)=0, so \eqref{eq:section_3_4} in this case becomes:

\noindent \begin{equation}
{\rm f}\left({\rm z}\right){\rm =2}{\rm z}^{{\rm 2}}{\rm \ }\int_{0}^{\infty}{\frac{\varphi{\rm (k)dk}}{{\rm k(}{\rm k}^{{\rm 2}}{\rm +}{\rm z}^{{\rm 2}}{\rm )}}}\label{eq:section_4_3}\end{equation}

\noindent with

\noindent \begin{equation}
\varphi\left({\rm k}\right){\rm =\ }\sum_{{{\rm k}\ge{\rm (2l-1)}\frac{\pi}{{\rm 2}}}_{{\rm \ }}}{\rm 1}\label{eq:section_4_4}\end{equation}
 We split \eqref{eq:section_4_4} up into a continuous part:

\noindent \[
{\varphi}_{{\rm 1}}\left({\rm k}\right){\rm =}\frac{{\rm k}}{\pi}\]

\noindent and a sawtooth type odd periodic function $\varphi_{2}(k)$
with period $\pi$ which over the interval $-\frac{\pi}{{\rm 2}}{\rm <}k<\frac{\pi}{{\rm 2}}$
equals:

\noindent \[
{\varphi}_{{\rm 2}}\left({\rm k}\right){\rm =-}\frac{{\rm k}}{\pi}\]

\noindent For $k>\frac{\pi}{{\rm 2}}$ $\varphi_{2}(k)$ just repeats
itself.

\noindent For our purposes we only need to consider $\varphi_{2}(k)$
in the range $0\le k<\infty$ .

\noindent \begin{equation}
\varphi\left({\rm k}\right){\rm =}\frac{{\rm k}}{\pi}{\rm +\ }{\varphi}_{{\rm 2}}\left({\rm k}\right)\label{eq:section_4_5}\end{equation}
 $\varphi_{2}(k)$ may be expanded in a Fourier series:

\noindent \begin{equation}
{\varphi}_{{\rm 2}}\left({\rm k}\right){\rm =}\sum_{n=1}^{\infty}{{\left({\rm -}{\rm 1}\right)}^{{\rm n}}\left.{{\rm \ }\frac{{\sin{\rm 2nk}\ }}{\pi{\rm n}}\ }\right.}\label{eq:section_4_6}\end{equation}
 From \eqref{eq:section_4_3}

\[
{\rm f}\left({\rm z}\right){\rm =2}{\rm z}^{{\rm 2}}\int_{0}^{\infty}{\frac{{\rm dk}}{{\rm k(}{\rm k}^{{\rm 2}}{\rm +}{\rm z}^{{\rm 2}}{\rm )}}\left[\frac{{\rm k}}{\pi}{\rm +}{\varphi}_{{\rm 2}}{\rm (k)}\right]}\]

\noindent The Dirichlet conditions apply to $\varphi_{2}(k)$, so
we are allowed to integrate term by term for $\varphi_{2}(k)$.

\noindent \begin{equation}
{\rm f}\left({\rm z}\right){\rm =}\frac{{\rm 2}{\rm z}^{{\rm 2}}}{\pi}\left[\int_{0}^{\infty}{\frac{{\rm dk}}{{\rm (}{\rm k}^{{\rm 2}}{\rm +}{\rm z}^{{\rm 2}}{\rm )}}}{\rm +\ }\sum_{{\rm n=1}}^{\infty}{\frac{{\rm (-1)}^{{\rm n}}}{{\rm n}}\int_{0}^{\infty}{\frac{{\sin\left({\rm 2nk}\right)\ }{\rm dk}}{{\rm k(}{\rm k}^{{\rm 2}}{\rm +}{\rm z}^{{\rm 2}}{\rm )}}}}\right]\label{eq:section_4_7}\end{equation}

\noindent The integral $\int_{0}^{\infty}{\frac{{\sin{\rm y}\ }dy}{y(y^{2}+a^{2})}}$
may be obtained via residue techniques and equals

\[
\frac{\pi}{2a^{2}}\ (1-e^{-a}).\]

\noindent Through appropriate change in variables:

\noindent \begin{equation}
\frac{{\rm 2}{\rm z}^{{\rm 2}}}{\pi}\sum_{{\rm n=1}}^{\infty}{\frac{{\rm (-1)}^{{\rm n}}}{{\rm n}}}\int_{0}^{\infty}{\frac{{\sin{\rm 2nk\ dk}\ }}{{\rm k(}{\rm k}^{{\rm 2}}{\rm +}{\rm z}^{{\rm 2}}{\rm )}}}{\rm =}\sum_{{\rm n=1}}^{\infty}{\frac{{\left({\rm -}{\rm 1}\right)}^{{\rm n}}}{{\rm n}}}\left({\rm 1-}{\rm e}^{{\rm -}{\rm 2nz}}\right){\rm =}\sum_{{\rm n=1}}^{\infty}{\frac{{\left({\rm -}{\rm 1}\right)}^{{\rm n}}}{{\rm n}}}{\rm +}\sum_{{\rm n=1}}^{\infty}{\frac{{\left({\rm -}{\rm 1}\right)}^{{\rm n+1}}}{{\rm n}}}{\rm e}^{{\rm -}{\rm 2nz}}\label{eq:section_4_8}\end{equation}

\noindent Both series are conditionally convergent and equal respectively
$-ln2$

\noindent and $ln(1+e^{-2z})$. \eqref{eq:section_4_7} then becomes:

\noindent $f(z)=z-ln2+ln(1+e^{-2z})$.

\noindent Through exponentiation and analytic continuation to the
entire z plane:

\noindent \[
{\rm F}\left({\rm z}\right){\rm =}\frac{{\rm e}^{{\rm z}}{\rm +}{\rm e}^{{\rm -}{\rm z}}}{{\rm 2}}\]

\section{\textit{The transform from $\varphi(k)$ to $f(z)$.}}

\noindent For a given $\varphi(k)$, the zeros of $f(z)$ are unaffected
by the assumed value of $F(0)$. Conveniently in this section we simply
put $F(0)=1$. \eqref{eq:section_3_4} then reduces to:

\noindent \begin{equation}
{\rm f}\left({\rm z}\right){\rm =2}{\rm z}^{{\rm 2}}{\rm \ }\int_{0}^{\infty}{\frac{\varphi{\rm (k)dk}}{{\rm k(}{\rm k}^{{\rm 2}}{\rm +}{\rm z}^{{\rm 2}}{\rm )}}}\label{eq:section_5_1}\end{equation}
 To the step function $\varphi(k)$ as defined by \eqref{eq:section_3_2}
corresponds a specific $f(z)$ as given by \eqref{eq:section_5_1}
and a specific $F(z)=e^{f(z)}$ as well.

\noindent Assuming $F(z)$ constitutes an axial product as defined
in this paper, $\varphi(k)$ may be recovered unambiguously from $F(z)$
by taking half of the integral

\[
\frac{{\rm 1}}{{\rm 2}\pi{\rm i}}\oint_{{\rm C}_{{\rm k}}}{\frac{{\rm F}^{{\rm '}}{\rm (z)}}{{\rm F(z)}}}{\rm dz}\]
 where $C_{k}$ is a contour dependent on $k$ in the sense that it
encircles the imaginary axis counterclockwise, crossing this axis
ones at each of the points $z=ik$ and $z=-ik$, e.g. a circle with
radius $k$ and centered at $z=0$.\\

\noindent We should in this context keep in mind that although $f(z)=lnF(z)$
obtained through \eqref{eq:section_5_1} is restricted to the semiplane
$x>0$, exponentiation to $F(z)$ and analytic continuation yields
an analytic function $F(z)$ over the entire plane.

\section{\textit{Transform table.}}

\noindent In table 1 are listed the results of the transform \eqref{eq:section_5_1}
for a number of different functions $\varphi(k)$. With the exception
of a possible initial discontinuity, all these functions are continuous
for the rest of the interval up to $\infty$. Although $\varphi(k)$
as defined in \eqref{eq:section_3_2} is non-decreasing, some of the
functions in table 1 are non-increasing or decreasing. By themselves
such $\varphi(k)$ do not make sense. Added to a non-decreasing $\varphi(k)$
however, they do make sense. Functions similar to $\varphi_{2}(k)$
in \eqref{eq:section_4_5}, which have discontinuities over the entire
range of $k$ are not included in the table. Although such contributions
to the total $\varphi(k)$ are essential in order to create zeros
for $F(z)$, for the asymptotic behaviour of $F(z)$ however they
just add a multiplicative constant. This is well illustrated by the
example of $Cosh\: z$ in section 4: For $\Re(z)\rightarrow\infty$,
$Cosh\: z\sim\frac{e^{z}}{2}$, whereas for the zeros of $Cosh\: z$,
the contribution of $\varphi_{2}(k)$ is essential. Without it, $\frac{e^{z}}{2}$
would not have any zeros. In addition it is essential for making $F(z)$
an even function. The $\varphi(k)$ listed in table 1 are more than
sufficient for our further needs. The function $u(k-a)$ in table
1 represents the unit step function, i.e. $u(k-a)=0$ for $k<a$,
while it equals 1 for $k\geq a$.\\
 \begin{tabular}{|c|c|c|l|}
\hline 
\multicolumn{1}{|||c}{Table 1} & \multicolumn{1}{c||||}{} &  & \tabularnewline
\hline 
 & $\varphi(k)$  &  & $f(z)$\tabularnewline
\hline 
1  & $u(k-a)$  & a>0  & ${ln\left(1+\frac{z^{2}}{a^{2}}\right)\ }=2ln{\rm z-2}ln{\rm a+O}\left({\rm z}^{{\rm -}{\rm 2}}\right)$\tabularnewline
\hline 
2  & $ku(k-a)$  & a>0  & $\pi{\rm z-2a+O}\left({\rm z}^{{\rm -}{\rm 2}}\right)$ \tabularnewline
\hline 
3  & \textit{$lnku(k-a)$}  & $\mid z\mid\geqslant1$; a>0  & ${(ln{\rm z})}^{2}{\rm +}\frac{{\pi}^{{\rm 2}}}{{\rm 12}}{\rm -}\frac{{\rm 1}}{{\rm 2}}\sum_{{\rm m=1}}^{\infty}{\frac{{\rm (-1)}^{{\rm m+1}}}{{\rm m}^{{\rm 2}}{\rm z}^{{\rm 2m}}}}$ \tabularnewline
\hline 
4  & $klnku(k-a)$  & $\mid z\mid\geqslant1$; a>0  & $\pi{\rm z}ln{\rm z-2a(}ln{\rm a-1)+\ O}\left({\rm z}^{{\rm -}{\rm 2}}\right)$ \tabularnewline
\hline 
5  & $\frac{ln{\rm k.u(k-a)}}{{\rm k}}$  & $\mid z\mid\geqslant1$; a>0  & ${\rm 2}\left[{\rm 1-}\sum_{{\rm m=1}}^{\infty}{\frac{{\rm (-1)}^{{\rm m+1}}}{{\rm (2m-1)}^{{\rm 2}}{\rm z}^{{\rm 2m}}}}\right]{\rm -}\frac{\pi ln{\rm z}}{{\rm z}}$ \tabularnewline
\hline 
6  & $k\sqrt{k}$  &  & $\pi\sqrt{{\rm 2}}{\rm z}\sqrt{{\rm z}}$ \tabularnewline
\hline 
7  & $k\sqrt{k}ln\: k$  &  & $\frac{\pi{\rm z}\sqrt{{\rm z}}}{\sqrt{{\rm 2}}}{\rm (2}ln{\rm z+}\pi{\rm )}$\tabularnewline
\hline 
8  & $\frac{u(k-a)}{k}$  & a>0  & $\frac{{\rm 2}}{{\rm a}}{\rm -}\frac{\pi}{{\rm z}}{\rm +}\frac{{\rm 2}}{{\rm z}}\left[{\rm arctg}\left(\frac{{\rm a}}{{\rm z}}\right)\right]$ \tabularnewline
\hline 
9  & $\frac{u(k-a)}{k^{2}}$  & a>0  & $\frac{1}{a^{2}}-\frac{1}{z^{2}}ln(1+\frac{z^{2}}{a^{2}}$) \tabularnewline
\hline
\end{tabular}

\section{\textit{Strip products}}

\noindent In strip products some zeros lie off-axis. The total product
equals $H(z)=F(z)G(z)$ in which $F(z)$ contains the axial zeros
and $G(z)$ the off-axis zeros.

\noindent Let $g(z)=ln\: G(z)$.

\noindent From \eqref{eq:section_2_5}:

\[
{\rm g}\left({\rm z}\right){\rm =}ln{\rm \ G}\left(0\right){\rm +}\sum_{{\rm m=1}}^{\infty}{{\rm (1+\ }\frac{{\rm 2}{\rm z}^{{\rm 2}}{\rm cos}{\beta}_{{\rm m}}}{{\rm q}_{{\rm m}}^{{\rm 2}}}{\rm +}\frac{{\rm z}^{{\rm 4}}}{{\rm q}_{{\rm m}}^{{\rm 4}}}{\rm )}}\]

\noindent For strip products we need to modify the definition of $ln\: G(z)$
slightly since zeros of $G(z)$ occur for $-\frac{1}{2}<x<\frac{1}{2}$.
We therefore consider the definition of $G(z)$ only for the halfplane
$x>\frac{1}{2}$. The formula used in section 2 for the logarithmic
function ${\rm g}\left({\rm z}\right){\rm =}{ln{\rm G(0)}}{\rm +}\int_{0}^{\infty}{\frac{{\rm G}^{{\rm '}}{\rm (w)}}{{\rm G(w)}}}dw$
is still valid if the path of integration starts at $z=0$, follows
the real axis up to $z=\frac{1}{2}$, and from there onwards to a
point z for which $x\geq\frac{1}{2}$

\noindent For the asymptotic behavior of $g(z)$ this restriction
poses no problem. Let's assume a quad of zeros with parameters $q_{m}$
and $\beta_{m}$. Let $C(q_{m})$ be the circle with radius $q_{m}$
and centered at $z=0$. Let $P$ be the intersection of $C(q_{m})$
with the vertical $z=\frac{1}{2}+it$, with $t>0$.

\noindent Let $\beta_{M}$ be the angle between the vertical $z=it$
and the radius from $z=0$ to $P$. Then:

\begin{equation}
{\cos{\beta}_{m}=\ \sqrt{1-\frac{1}{4q_{m}^{2}}}}\label{eq:section_7_1}\end{equation}

\noindent The first zero in the quad (see (2.4)), i.e. $iq_{m}e^{-i{\beta}_{m}}$
will be located on $C(q_{m})$ and the associated angle $\beta_{m}$
will be less than $\beta_{M}$ . Therefore:

\noindent \begin{equation}
cos\:\beta{}_{{\rm m}}{\rm >}cos\:\beta_{M}\label{eq:section_7_2}\end{equation}

\noindent \eqref{eq:section_7_1} implies \begin{equation}
cos\:\beta_{M}>1-\frac{1}{8q_{m}^{2}}\label{eq:section_7_3}\end{equation}

\noindent \eqref{eq:section_7_2} and \eqref{eq:section_7_3} imply:

\[
{{\rm 1}>cos{\beta}_{{\rm m}}\ }{\rm >}1-\frac{1}{8q_{m}^{2}}\]
 or: \begin{equation}
{\cos{\beta}_{{\rm m}}\ }{\rm =}1-\frac{{\theta}_{m}}{8q_{m}^{2}}\quad with\;0<\theta_{m}<1\label{eq:section_7_4}\end{equation}

\noindent \eqref{eq:section_7_4}implies:\begin{equation}
ln\left({\rm 1+}\frac{{\rm 2}{\rm z}^{{\rm 2}}{\rm cos}{\beta}_{{\rm m}}}{{\rm q}_{{\rm m}}^{{\rm 2}}}{\rm +}\frac{{\rm z}^{{\rm 4}}}{{\rm q}_{{\rm m}}^{{\rm 4}}}\right){\rm =2}ln\left({\rm 1+}\frac{{\rm z}^{{\rm 2}}}{{\rm q}_{{\rm m}}^{{\rm 2}}}\right){\rm +ln}\left[{\rm 1-}\frac{{\rm z}^{{\rm 2}}{\theta}_{{\rm m}}}{{{\rm 4}{\rm q}_{{\rm m}}^{{\rm 4}}{\rm (1+}\frac{{\rm z}^{{\rm 2}}}{{\rm q}_{{\rm m}}^{{\rm 2}}}{\rm )}}^{{\rm 2}}}\right]\label{eq:Section_7_5}\end{equation}

\noindent \[
\frac{{\rm z}^{{\rm 2}}{\theta}_{{\rm m}}}{{{\rm 4}{\rm q}_{{\rm m}}^{{\rm 4}}\left({\rm 1+}\frac{{\rm z}^{{\rm 2}}}{{\rm q}_{{\rm m}}^{{\rm 2}}}\right)}^{{\rm 2}}}{\rm =}\frac{{\theta}_{{\rm m}}}{{{\rm 4}{\rm q}_{{\rm m}}^{{\rm 4}}\left(\frac{{\rm 1}}{{\rm z}^{{\rm 2}}}{\rm +}\frac{{\rm 2}}{{\rm q}_{{\rm m}}^{{\rm 2}}}{\rm +}\frac{{\rm z}^{{\rm 2}}}{{\rm q}_{{\rm m}}^{{\rm 4}}}\right)}}\]

\noindent For $\Re(z)\to\infty$, we assume \begin{equation}
\Re(z^{2})>0\: i.e.\left|{\arg\left(z\right)}\right|<\frac{\pi}{4}\label{eq:section_7_6}\end{equation}

\noindent \eqref{eq:section_7_6} implies $ln\left[1-\frac{z^{2}{\theta}_{m}}{{4q_{m}^{4}(1+\frac{z^{2}}{q_{m}^{2}})}^{2}}\right]$
will go to zero for $\Re(z)\to\infty$.

\noindent Therefore all finite sums of such terms will equally go
to zero as $O(z^{-2})$.

\noindent Let's now take the case of an infinite sum of $q_{m}'s$.

\noindent Under constraint \eqref{eq:section_7_6}:$\left|\frac{{\theta}_{m}}{4\left[\frac{q_{m}^{4}}{z^{2}}+2q_{m}^{2}+z^{2}\right]}\right|<\ \frac{{\theta}_{m}}{8q_{m}^{2}}$

\noindent It is not difficult to prove that in this case and for $q_{m}>\frac{1}{2}$

\noindent $\sum_{m=1}^{\infty}{\left|ln\left[1-\frac{z^{2}{\theta}_{m}}{{4q_{m}^{4}(1+\frac{z^{2}}{q_{m}^{2}})}^{2}}\right]\right|{\rm <}K}$
with K some positive constant.

\noindent For a strip product \textit{$ln\: H(z)=ln\: H(0)+f(z)+g(z)$}

\noindent with:

\[
{\rm f}\left({\rm z}\right)=\ \sum_{{\rm l=1}}^{\infty}{ln\:{\rm (1+}\frac{{\rm z}^{{\rm 2}}}{{\rm k}_{{\rm l}}^{{\rm 2}}}})\]

\[
{\rm g}\left({\rm z}\right){\rm =}\sum_{{\rm m=1}}^{\infty}{{\rm 2}ln\:{\rm (1+}}\frac{{\rm z}^{{\rm 2}}}{{\rm q}_{{\rm m}}^{{\rm 2}}}{\rm )+\ }\sum_{{\rm m=1}}^{\infty}{ln\left[{\rm 1-}\frac{{\rm z}^{{\rm 2}}{\theta}_{{\rm m}}}{{{\rm 4}{\rm q}_{{\rm m}}^{{\rm 4}}{\rm (1+}\frac{{\rm z}^{{\rm 2}}}{{\rm q}_{{\rm m}}^{{\rm 2}}}{\rm )}}^{{\rm 2}}}\right]}\]
 As discussed above, the contribution of the second term of $g(z)$
will be $O(z^{-2})$ for a finite number of terms and at most $O(1)$
for an infinite number of terms.

\noindent The first term of $g(z)$ is simply the equivalent of a
double axial product term, i.e. each $q_{m}$ may be treated as a
jump of size 2 in the $\varphi(k)$ function. In fact, the same argument
applies to multiple roots. Multiple roots will simply cause larger
jumps in $\varphi(k)$.

\noindent As an example of multiple roots we consider an exercise
similar to the one for $Cosh\: z$ (section 4). Here we consider $\varphi\left({\rm k}\right){\rm =}\frac{{\rm k}}{\pi}{\rm +}{\varphi}_{{\rm 2}}\left({\rm k}\right)$
where for $-\frac{N\pi}{2}<k<\frac{N\pi}{2}$, ${\varphi}_{{\rm 2}}\left({\rm k}\right){\rm =-}\frac{{\rm k}}{\pi}$
and for $\frac{N\pi}{2}<k$, ${\varphi}_{{\rm 2}}\left({\rm k}\right)$
repeats itself with period $N\pi$. $N$ being a positive integer.
The resulting $\varphi(k)$ has on the average over a given interval
$0<k<K$, the same number of zeros as for $Cosh\: z$, but now we
bunch $N$ zeros together, i.e. the size of the jumps equals $N$
instead of 1. Each zero has multiplicity $N$.

\noindent A similar calculation then yields:

\[
{\rm f}\left({\rm z}\right){\rm =z-N}ln{\rm 2+N}ln{\rm (1+}{\rm e}^{{\rm -}{\rm 2}\frac{{\rm z}}{{\rm N}}}{\rm )}\]

\[
{\rm F}\left({\rm z}\right){\rm =\ }{\rm e}^{{\rm f(z)}}{\rm =}{\left[{{\rm Cosh}\left(\frac{{\rm z}}{{\rm N}}\right)\ }\right]}^{{\rm N}}\]
 For $z\rightarrow\infty$: \[
\frac{{\left[{\rm Cosh}\left(\frac{{\rm z}}{{\rm N}}\right)\right]}^{{\rm N}}}{{\rm Cosh\ z}}{\rm \sim}\frac{{\rm 1}}{{\rm 2}^{{\rm N-1}}}\]

\noindent So for this example, bunching the zeros together leads to
an asymptotic reduction of $F(z)$ equal to $\frac{1}{2^{N-1}}$.

\section{\textit{Asymptotic behavior of ln $\xi(z)$.}}

\noindent Note: It would be misleading to write $\xi$(z+$\frac{1}{2}$)
for $\xi(s)$ after substitution $s=z+\frac{1}{2}$. Indeed $\xi(z)$
is an even function of $z$, not of $z+\frac{1}{2}$.

\noindent \[
\xi\left({\rm s}\right){\rm =\ }\frac{{\rm 1}}{{\rm 2}}\Gamma\left(\frac{{\rm s}}{{\rm 2}}\right){\pi}^{{\rm -}\frac{{\rm s}}{{\rm 2}}}{\rm \ s(s-1)}\zeta{\rm (s)}\]

\noindent For \textit{ln} $\Gamma(\frac{s}{2})$ we use Stirling's
formula (see Lang, p. 423):

\begin{equation}
ln\Gamma\left({\rm A}\right){\rm =}\left({\rm A-}\frac{{\rm 1}}{{\rm 2}}\right)ln{\rm A-A+}\frac{{\rm 1}}{{\rm 2}}{ln\left({\rm 2}\pi\right)\ }{\rm -}{\rm w}\left({\rm A}\right)\label{eq:section_8_1}\end{equation}
 For A we will need to plug in ${\rm A=}\frac{{\rm s}}{{\rm 2}}{\rm =}\frac{{\rm z}}{{\rm 2}}{\rm +}\frac{{\rm 1}}{{\rm 4}}$

\noindent The function $w(A)$ equals:

\begin{equation}
{\rm w}\left({\rm A}\right){\rm =}\int_{0}^{\infty}{\frac{{\rm P(}\lambda{\rm )d}\lambda}{{\rm A+}\lambda}}\label{eq:section_8_2}\end{equation}
 in which $P\left(\lambda\right)$is a sawtooth type odd periodic
function with period 1. In the interval $0<\lambda<1$ $P(\lambda)$
equals $-\frac{1}{2}+\lambda$. After that it repeats itself.

\noindent $w(A)$is analytic in the complex $A$ plane from which
the origin and the negative real axis is removed.

\noindent $w(A)$may be written in various forms:

\noindent \[
{\rm w}\left({\rm A}\right){\rm =-}\sum_{{\rm r=1}}^{\infty}{\sum_{{\rm m=2}}^{\infty}{\left[\frac{{\rm 1}}{{\rm m+1}}{\rm -}\frac{{\rm 1}}{{\rm 2m}}\right]}}\frac{{\rm 1}}{{\rm (A+r)}^{{\rm m}}}\]

\noindent or more condensed (see Freitag-Busam, p204):

\noindent \[
{\rm w}\left({\rm A}\right){\rm =-}\sum_{{\rm k=0}}^{\infty}{\left[\left({\rm A+k+}\frac{{\rm 1}}{{\rm 2}}\right){\ln\left({\rm 1+}\frac{{\rm 1}}{{\rm A+k}}\right)\ }{\rm -}{\rm 1}\right]}\]

\noindent or using the generalized or Hurwitz zeta function,

\noindent \[
{\zeta}_{{\rm g}}{\rm (2m,A+}\frac{{\rm 1}}{{\rm 2}}{\rm )=}\sum_{{\rm k=0}}^{\infty}{{{\rm (A+}\frac{{\rm 1}}{{\rm 2}}{\rm +k)}}^{{\rm -}{\rm 2m}}}\]

\[
{\rm w}\left({\rm A}\right){\rm =}\sum_{{\rm m=1}}^{\infty}{\frac{{\zeta}_{{\rm g}}{\rm (2m,A+}\frac{{\rm 1}}{{\rm 2}}{\rm )}}{{\rm (2m+1)}{\rm 2}^{{\rm 2m}}}}\]
 It is easy to demonstrate that for $\Re(z)>8\lyxmathsym{\textonehalf}$:

\noindent \begin{equation}
\left|{\rm w}\left[\frac{{\rm z}}{{\rm 2}}{\rm +}\frac{{\rm 1}}{{\rm 4}}\right]\right|{\rm <}\frac{{\rm 1}}{{\rm 8x}}{\rm \ \ \ \ \ x=\Re(z)}\label{eq:section_8_3}\end{equation}

\noindent \eqref{eq:section_8_1}yields:\begin{equation}
{ln\Gamma\left(\frac{{\rm z}}{{\rm 2}}{\rm +}\frac{{\rm 1}}{{\rm 4}}\right)\ }{\rm =}\frac{{\rm z}}{{\rm 2}}{ln\left(\frac{{\rm z}}{{\rm 2}}\right)\ }{\rm -}\frac{{\rm 1}}{{\rm 4}}{ln\left(\frac{{\rm z}}{{\rm 2}}\right)\ }{\rm +}\frac{{\rm z}}{{\rm 2}}{ln\left({\rm 1+}\frac{{\rm 1}}{{\rm 2z}}\right)\ }{\rm -}\frac{{\rm 1}}{{\rm 4}}{ln\left({\rm 1+}\frac{{\rm 1}}{{\rm 2z}}\right){\rm -}\frac{{\rm z}}{{\rm 2}}{\rm -}\frac{{\rm 1}}{{\rm 4}}{\rm +}\frac{{\rm 1}}{{\rm 2}}{ln\left({\rm 2}\pi\right){\rm -}{\rm w}\left(\frac{{\rm z}}{{\rm 2}}{\rm +}\frac{{\rm 1}}{{\rm 4}}\right)}}\label{eq:section_8_4}\end{equation}

\noindent Further:

\begin{equation}
{ln\xi\left({\rm z}\right){\rm =-}ln{\rm 2+}ln\Gamma{\rm (}\ }\frac{{\rm z}}{{\rm 2}}{\rm +}\frac{{\rm 1}}{{\rm 4}}{\rm )-}\left(\frac{{\rm z}}{{\rm 2}}{\rm +}\frac{{\rm 1}}{{\rm 4}}\right)ln\pi{\rm +ln}\left({\rm z}^{{\rm 2}}{\rm -}\frac{{\rm 1}}{{\rm 4}}\right){\rm +}ln\zeta\left({\rm z+}\frac{{\rm 1}}{{\rm 2}}\right)\label{eq:section_8_5}\end{equation}
 Due to the presence of $ln\zeta\left(z+\frac{1}{2}\right)$and ${\ln\left(z^{2}-\frac{1}{4}\right)\ }$
we define ln$\xi$(z) only for the semiplane $x>\frac{1}{2}$, similar
to what is done in section 7, but excluding $z=\frac{1}{2}$ itself.

\noindent Considering that for $x>10$ we have $|ln\zeta{\rm (z+1/2)}|$$<$$\frac{20}{19x}$
and putting some bounds on series expansions of certain logarithmic
terms, we obtain:

\begin{equation}
ln\xi\left({\rm z}\right){\rm =}\frac{{\rm z}}{{\rm 2}}{ln\left(\frac{{\rm z}}{{\rm 2}\pi}\right)\ }{\rm -}\frac{{\rm z}}{{\rm 2}}{\rm +}\frac{{\rm 7}}{{\rm 4}}ln{\rm z+}\frac{{\rm 1}}{{\rm 4}}{ln\left(\frac{\pi}{{\rm 2}}\right)\ }{\rm +}\frac{{\rm 2}\vartheta_{{\rm 1}}}{{\rm z}}\label{eq:section_8_6}\end{equation}
 with $\left|\vartheta_{1}\right|<1,\Re(z)>10$

\noindent At this point one should not get the impression that $\zeta(s)$
is not significant since its contribution to the asymptotic form of
ln $\xi(z)$ disappears in a term $\frac{\vartheta_{1}}{z}$. Indeed,
$\zeta(s)$ is an essential component in \eqref{eq:section_1_1} to
create all the zeros of $\xi(s)$ and making it an even function.

\section{\textit{Distribution of the roots for $\xi(z)$.}}

\noindent The intention now is to find an approximate formula for
the number of roots up to a certain k value, based on \eqref{eq:section_8_6}.
To keep things simple, in what follows we will consider real values
of z. This will not constrain our conclusions.

\noindent \eqref{eq:section_8_6} contains three terms which keep
increasing without bound in absolute value with increasing z:

\noindent \[
T_{1}=\frac{z}{2}{\ln\left(\frac{z}{2\pi}\right)\ };T_{2}=-\frac{z}{2};T_{3}=\frac{7}{4}{\ln\: z}\]

\noindent In addition there is a constant and a vanishing term for
$z\rightarrow\infty$:

\noindent \[
\frac{{\rm 1}}{{\rm 4}}{ln\left(\frac{\pi}{{\rm 2}}\right)}{\rm +}\frac{{\rm 2}\vartheta_{{\rm 1}}}{{\rm z}}\]

\noindent If $\varphi(k)$ is the actual step function representing
the distribution of the roots of $\xi(z)$, we can split $\varphi(k)$
up into a continuous curve $\phi(k)$ and a discontinuous sawtooth
type function $\Omega(k)$:

\begin{equation}
\varphi\left({\rm k}\right){\rm =}\phi\left({\rm k}\right){\rm +}\Omega\left({\rm k}\right)\label{eq:section_9_1}\end{equation}
 $\frac{{\rm d}\Omega{\rm (k)}}{{\rm dk}}{\rm =-}\frac{{\rm d}\phi{\rm (k)}}{{\rm dk}}$
for almost all $k$ except for $k=k_{l}$ values corresponding to
roots of $\xi(z)$. At those $k_{l}$, jumps will occur in $\Omega(k)$
equal to the multiplicity of the root. We will first concentrate on
the terms $T_{1},T_{2},T_{3}$.

\noindent The idea is to identify for each of these terms a corresponding
term in table 1. We are of course allowed to multiply the $\varphi(k)$
in table 1 with appropriate coefficients in order to obtain the required
correspondence with \eqref{eq:section_8_6}.

\noindent To:

\[
\frac{{\rm k}}{{\rm 2}\pi}ln{\rm k.u}\left({\rm k-a}\right){\rm \ \ \ \ corresponds\ \ \ \ \ \ \ \ \ \ }\frac{{\rm z}}{{\rm 2}}ln{\rm z-}\frac{{\rm a}}{\pi}\left(ln{\rm a-1}\right){\rm +O(}{\rm z}^{{\rm -}{\rm 2}}{\rm )}\]

\[
{\rm -}\frac{{\rm k}}{{\rm 2}\pi}ln\left({\rm 2}\pi\right){\rm .u}\left({\rm k-a}\right){\rm \ \ \ \ \ \ corresponds\ \ \ \ \ -}\frac{{\rm z}}{{\rm 2}}ln{\rm 2}\pi{\rm +}\frac{{\rm a}}{\pi}ln{\rm 2}\pi{\rm +O(}{\rm z}^{{\rm -}{\rm 2}}{\rm )}\]

\[
{\rm -}\frac{{\rm k}}{{\rm 2}\pi}{\rm u}\left({\rm k-a}\right){\rm \ \ \ \ \ \ \ \ \ \ \ \ \ \ \ \ \ \ \ \ \ \ \ \ \ \ \ \ corresponds\ \ \ \ \ \ \ \ \ \ \ \ \ \ -}\frac{{\rm z}}{{\rm 2}}{\rm +}\frac{{\rm a}}{\pi}{\rm +O(}{\rm z}^{{\rm -}{\rm 2}}{\rm )}\]

\[
\frac{{\rm 7}}{{\rm 8}}{\rm u}\left({\rm k-a}\right){\rm \ \ \ \ \ \ \ \ \ \ \ corresponds\ \ \ \ \ \ \ \ \ \ }\frac{{\rm 7}}{{\rm 4}}ln{\rm z-}\frac{{\rm 7}}{{\rm 4}}ln{\rm a+O(}{\rm z}^{{\rm -}{\rm 2}}{\rm )}\]

In total to:

\begin{eqnarray*}
\left[\frac{{\rm k}}{{\rm 2}\pi}{ln\left(\frac{{\rm k}}{{\rm 2}\pi}\right)\ }{\rm -}\frac{{\rm k}}{{\rm 2}\pi}{\rm +}\frac{{\rm 7}}{{\rm 8}}\right]{\rm .u}\left({\rm k-a}\right) &  & corresponds\\
 &  & \frac{{\rm z}}{{\rm 2}}{ln\left(\frac{{\rm z}}{{\rm 2}\pi}\right)\ }{\rm -}\frac{{\rm z}}{{\rm 2}}{\rm +}\frac{{\rm 7}}{{\rm 4}}ln{\rm z+}\frac{{\rm a}}{\pi}\left[{\rm 2-}ln{\rm a+}ln{\rm 2}\pi\right]{\rm -}\frac{{\rm 7}}{{\rm 4}}ln{\rm a+O(}{\rm z}^{{\rm -}{\rm 2}}{\rm )}\end{eqnarray*}

\noindent For $k>a$ we may drop $u(k-a)$. To shorten things we call:

\noindent \begin{equation}
\phi\left({\rm k}\right){\rm =}\frac{{\rm k}}{{\rm 2}\pi}{ln\left(\frac{{\rm k}}{{\rm 2}\pi}\right)}{\rm -}\frac{{\rm k}}{{\rm 2}\pi}{\rm +}\frac{{\rm 7}}{{\rm 8}}\label{eq:section_9_2}\end{equation}

\begin{equation}
{\rm T}_{{\rm 4}}{\rm (z)=}\frac{{\rm z}}{{\rm 2}}{ln\left(\frac{{\rm z}}{{\rm 2}\pi}\right)\ }{\rm -}\frac{{\rm z}}{{\rm 2}}{\rm +}\frac{{\rm 7}}{{\rm 4}}ln\:{\rm z}\label{eq:section_9_3}\end{equation}

\begin{equation}
{\rm T}_{{\rm 5}}\left({\rm z}\right){\rm =}{\rm T}_{{\rm 4}}\left({\rm z}\right){\rm +}\frac{{\rm a}}{\pi}\left[{\rm 2-}ln{\rm a+}ln{\rm 2}\pi\right]{\rm -}\frac{{\rm 7}}{{\rm 4}}ln{\rm a+\ O(}{\rm z}^{{\rm -}{\rm 2}}{\rm )}\label{eq:section_9_4}\end{equation}

\noindent We select the value of a such that we start with $\phi$(a)=0.
This happens for $k=a=9.6769$...

\noindent $T_{5}\left(z\right)$ then becomes: \begin{equation}
T_{5}\left(z\right)=T_{4}\left(z\right)+0,8582\dots+O\left(z^{-2}\right)\label{eq:section_9_5}\end{equation}

\noindent \begin{flushleft}
We now compare $T_{5}\left(z\right)$ with ln$\xi$(z), keeping in
mind that 
\par\end{flushleft}

\noindent \begin{flushleft}
\begin{equation}
ln\xi\left({\rm z}\right){\rm =}ln\xi\left(0\right){\rm +}ln\prod_{{\rm l=1}}^{\infty}{\left({\rm 1+}\frac{{\rm z}^{{\rm 2}}}{{\rm k}_{{\rm l}}^{{\rm 2}}}\right)}\label{eq:section_9_6}\end{equation}

\par\end{flushleft}

\noindent ln $\xi(0)$ may be computed from \eqref{eq:section_1_1}:
$ln\xi(0)=-0,69892$...

\noindent From \eqref{eq:section_8_6} and \eqref{eq:section_9_6}:

\noindent \[
ln\prod_{{\rm l=1}}^{\infty}{\left({\rm 1+}\frac{{\rm z}^{{\rm 2}}}{{\rm k}_{{\rm l}}^{{\rm 2}}}\right){\rm -}}{\rm T}_{{\rm 5}}\left({\rm z}\right){\rm =\ }\frac{{\rm 1}}{{\rm 4}}{ln\left(\frac{\pi}{{\rm 2}}\right)\ }{\rm -}ln\xi\left(0\right){\rm -}{\rm 0,8582+\ O}\left({\rm z}^{{\rm -}{\rm 1}}\right)\]

\begin{equation}
ln\prod_{{\rm l=1}}^{\infty}{\left({\rm 1+}\frac{{\rm z}^{{\rm 2}}}{{\rm k}_{{\rm l}}^{{\rm 2}}}\right){\rm -}}{\rm T}_{{\rm 5}}\left({\rm z}\right){\rm =0,0464+O}\left({\rm z}^{{\rm -}{\rm 1}}\right)\label{eq:section_9_7}\end{equation}
 The remaining constant 0,0464 in \eqref{eq:section_9_7} must be
due to $\Omega(k)$, i.e.

\noindent ${{\lim}_{x\to\infty}2\int_{a}^{\infty}{\frac{\Omega\left(k\right)dk}{k(1+\frac{k^{2}}{z^{2}})}}}$
should equal 0,0464.

\section{\textit{The discontinuous function $\Omega(k)$.}}

\noindent In line with earlier definitions (see(3.2),\eqref{eq:section_9_1})
we may put $\Omega(k)=0$ for $0\leqslant k\leqslant a$ in which
a is the point where $\phi(k)$ starts off from zero. According to
the foregoing, $a<k_{1}$ where $k_{1}$ corresponds to the first
zero of $\xi$(z).

\noindent Let in general the transform $\mathfrak{T}$ from the k
domain to the z domain be understood as:

\begin{equation}
\mathfrak{T}\left[\varphi{\rm (k)}\right]{\rm =2}{\rm z}^{{\rm 2}}\int_{{\rm a}}^{\infty}{\frac{\Omega\left({\rm k}\right){\rm dk}}{{\rm k(}{\rm k}^{{\rm 2}}{\rm +}{\rm z}^{{\rm 2}}{\rm )}}}\label{eq:section_10_1}\end{equation}

\noindent Let:

\begin{equation}
{{\lim}_{{\rm T}\to\infty}\frac{{\rm 1}}{{\rm T}}\int_{{\rm a}}^{{\rm T}}{\Omega\left({\rm k}\right){\rm dk=}{\Omega}_{0}}}\label{eq:section_10_2}\end{equation}
 Let also:

\noindent \begin{equation}
\Omega\left({\rm k}\right){\rm =}{\Omega}_{{\rm 1}}\left({\rm k}\right){\rm +}{\Omega}_{0}\left({\rm k}\right){\rm \ with\ }{\Omega}_{0}\left({\rm k}\right){\rm =}{\Omega}_{0}{\rm u}\left({\rm k-a}\right)\label{eq:section_10_3}\end{equation}

\begin{equation}
\mathfrak{T}\left[\Omega\left({\rm k}\right)\right]{\rm =}\mathfrak{T}\left[{\Omega}_{{\rm 1}}\left({\rm k}\right)\right]{\rm +}\mathfrak{T}\left[{\Omega}_{0}{\rm u}\left({\rm k-a}\right)\right]{\rm =\mathfrak{T}\left[{\Omega}_{{\rm 1}}\left({\rm k}\right)\right]{\rm +}{\Omega}_{0}{\ln\left({\rm 1+}\frac{{\rm z}^{{\rm 2}}}{{\rm a}^{{\rm 2}}}\right)\ }}\label{eq:section_10_4}\end{equation}
 The existence of the logarithmic term in \eqref{eq:section_10_4}
as part of ln$\xi$(z) is in contradiction with \eqref{eq:section_9_7}
where no logarithmic term is left. So $\Omega_{0}$ must equal zero.

\noindent Therefore by definition:

\begin{equation}
{{\lim}_{{\rm T}\to\infty}{\int_{{\rm a}}^{{\rm T}}{{\Omega}_{{\rm 1}}{\rm (k)dk=0}}}}\label{eq:section_10_5}\end{equation}
 If ${\Omega}_{{\rm 1}}\left({\rm k}\right)$were a smooth function,
it could after one or more zero-line crossings remain permanently
above or below this line and approach it asymptotically and still
meet \eqref{eq:section_10_5}.

\noindent But ${\Omega}_{{\rm 1}}\left({\rm k}\right)$is a discontinuous
function with jumps of at least height one. Therefore ${\Omega}_{{\rm 1}}\left({\rm k}\right)$
will cross the zero-line an infinite number of times. If we now look
back at the situation in terms of $\varphi(k)$ and $\phi(k)$ (see
\eqref{eq:section_9_1}), $\varphi(k)$ will cross $\phi(k)$ alternatively
vertically and horizontally.

\noindent At the horizontal crossings the number of zeros of $\xi(z)$
from $k=0$ to the point of crossing $k_{a}$ is equal to $\phi(k_{a})$.
If two consecutive horizontal crossings occur at $k_{a}$ and $k_{b}$,
the number of zeros in between will equal:

\noindent \[
\phi\left({\rm k}_{{\rm b}}\right){\rm -}\phi{\rm (}{\rm k}_{{\rm a}}{\rm )\sim}\frac{\triangle{\rm (a,b)}}{{\rm 2}\pi}{\rm ln}\left[\frac{{\rm k}_{{\rm a}}{\rm +}{\rm k}_{{\rm b}}}{{\rm 4}\pi}\right]\]
 with $\triangle\left({\rm a,b}\right){\rm =}{\rm k}_{{\rm b}}{\rm -}{\rm k}_{{\rm a}}$
.

\noindent At the horizontal crossings the total number of zeros from
0 up to T will equal:

\begin{equation}
N\left(T\right)=\frac{T}{2\pi}{\ln\left(\frac{T}{2\pi}\right)}-\frac{T}{2\pi}+\frac{7}{8}\label{eq:section_10_6}\end{equation}
 In \eqref{eq:section_10_6} the symbol T has been used instead of
k, as customary in analytic number theory.

\noindent The reader is invited to go through an exercise in graphics.
Draw the following curves in a diagram: 
\begin{enumerate}
\item $\phi(k)$ starting with $k=9.67$... and up to e.g. $k=55$. 
\item $\varphi_{1}(k)$: The $\varphi(k)$ corresponding to $\xi(z)$ based
on the actual zeros available in the literature (e.g. Edwards, p.
96) up to e.g. the tenth root. $\varphi_{1}(k)$ is a step function
with vertical jumps at the zeros of $\xi(z)$. The horizontal portions
occur for vertical values equal to n=0,1,2,... 
\item $\varphi_{2}(k)$: A step function similar to $\varphi_{1}(k)$ with
horizontal portions at the same level n=0,1,2,... as $\varphi_{1}(k)$,
while the vertical jumps occur at those horizontal values $k$ for
which $\phi(k)$ reaches a value $n+\frac{1}{2}$. 
\end{enumerate}
\noindent It is remarkable how well $\varphi_{2}(k)$ approximates
$\varphi_{1}(k)$. 

\noindent Renaat Van Malderen

\noindent Address: Maxlaan 21, 2640 Mortsel, Belgium

\noindent The author can be reached by email: hans.vanmalderen@skynet.be 
\end{document}